\documentclass{article}

\usepackage{amssymb,amsmath,theorem,euscript}

\input{epsf}

\def\sm{\smallskip}

\newcounter{punct}

\def\punct{\refstepcounter{punct}{\arabic{punct}. }}

\newtheorem{theorem}{Theorem}%{Теорема}
\newtheorem{lemma}[theorem]{Lemma}%{Лемма}

\begin{document}

 \def\ov{\overline}
\def\wt{\widetilde}
\def\wh{\widehat}
 \newcommand{\rk}{\mathop {\mathrm {rk}}\nolimits}
\newcommand{\Aut}{\mathop {\mathrm {Aut}}\nolimits}
\newcommand{\Out}{\mathop {\mathrm {Out}}\nolimits}
\renewcommand{\Re}{\mathop {\mathrm {Re}}\nolimits}
\renewcommand{\Im}{\mathop {\mathrm {Im}}\nolimits}
 \newcommand{\tr}{\mathop {\mathrm {tr}}\nolimits}
  \newcommand{\Hom}{\mathop {\mathrm {Hom}}\nolimits}
   \newcommand{\diag}{\mathop {\mathrm {diag}}\nolimits}
   \newcommand{\supp}{\mathop {\mathrm {supp}}\nolimits}
 \newcommand{\im}{\mathop {\mathrm {im}}\nolimits}

\def\Br{\mathrm {Br}}

 %%%1. ClASSICAL GROUPS
\def\SL{\mathrm {SL}}
\def\SU{\mathrm {SU}}
\def\GL{\mathrm {GL}}
\def\U{\mathrm U}
\def\OO{\mathrm O}
 \def\Sp{\mathrm {Sp}}
  \def\Ad{\mathrm {Ad}}
 \def\SO{\mathrm {SO}}
\def\SOS{\mathrm {SO}^*}
 \def\Diff{\mathrm{Diff}}
 \def\Vect{\mathfrak{Vect}}
\def\PGL{\mathrm {PGL}}
\def\PU{\mathrm {PU}}
\def\PSL{\mathrm {PSL}}
\def\Symp{\mathrm{Symp}}
\def\End{\mathrm{End}}
\def\Mor{\mathrm{Mor}}
\def\Aut{\mathrm{Aut}}
 \def\PB{\mathrm{PB}}
\def\Fl{\mathrm {Fl}}
\def\Symm{\mathrm {Symm}}

 \def\cA{\mathcal A}
\def\cB{\mathcal B}
\def\cC{\mathcal C}
\def\cD{\mathcal D}
\def\cE{\mathcal E}
\def\cF{\mathcal F}
\def\cG{\mathcal G}
\def\cH{\mathcal H}
\def\cJ{\mathcal J}
\def\cI{\mathcal I}
\def\cK{\mathcal K}
 \def\cL{\mathcal L}
\def\cM{\mathcal M}
\def\cN{\mathcal N}
 \def\cO{\mathcal O}
\def\cP{\mathcal P}
\def\cQ{\mathcal Q}
\def\cR{\mathcal R}
\def\cS{\mathcal S}
\def\cT{\mathcal T}
\def\cU{\mathcal U}
\def\cV{\mathcal V}
 \def\cW{\mathcal W}
\def\cX{\mathcal X}
 \def\cY{\mathcal Y}
 \def\cZ{\mathcal Z}
%%% END MATHCAL %%%%%%%%%%%%%%%%%%%%%%%%%%%%%%%%% %%%%%%%%%%%%%%%%%%%%%%%%%%%%%%%% %%%
\def\0{{\ov 0}}
 \def\1{{\ov 1}}
 %%%%%%%%%%%%%%%%%%%%%%%%%%%% %%%%%%%%%%%%%%%%%%%%%%%%%%%%%%%%%%% %%% BEGIN GOTIC
 \def\frA{\mathfrak A}
 \def\frB{\mathfrak B}
\def\frC{\mathfrak C}
\def\frD{\mathfrak D}
\def\frE{\mathfrak E}
\def\frF{\mathfrak F}
\def\frG{\mathfrak G}
\def\frH{\mathfrak H}
\def\frI{\mathfrak I}
 \def\frJ{\mathfrak J}
 \def\frK{\mathfrak K}
 \def\frL{\mathfrak L}
\def\frM{\mathfrak M}
 \def\frN{\mathfrak N} \def\frO{\mathfrak O} \def\frP{\mathfrak P} \def\frQ{\mathfrak Q} \def\frR{\mathfrak R}
 \def\frS{\mathfrak S} \def\frT{\mathfrak T} \def\frU{\mathfrak U} \def\frV{\mathfrak V} \def\frW{\mathfrak W}
 \def\frX{\mathfrak X} \def\frY{\mathfrak Y} \def\frZ{\mathfrak Z} \def\fra{\mathfrak a} \def\frb{\mathfrak b}
 \def\frc{\mathfrak c} \def\frd{\mathfrak d} \def\fre{\mathfrak e} \def\frf{\mathfrak f} \def\frg{\mathfrak g}
 \def\frh{\mathfrak h} \def\fri{\mathfrak i} \def\frj{\mathfrak j} \def\frk{\mathfrak k} \def\frl{\mathfrak l}
 \def\frm{\mathfrak m} \def\frn{\mathfrak n} \def\fro{\mathfrak o} \def\frp{\mathfrak p} \def\frq{\mathfrak q}
 \def\frr{\mathfrak r} \def\frs{\mathfrak s} \def\frt{\mathfrak t} \def\fru{\mathfrak u} \def\frv{\mathfrak v}
 \def\frw{\mathfrak w} \def\frx{\mathfrak x} \def\fry{\mathfrak y} \def\frz{\mathfrak z} \def\frsp{\mathfrak{sp}}
 %% This is Lie algebra %%% END GOTIC
%%%%%%%%%%%%%%%%%%%%%%%%%%%%%%%% %%%%%%%%%%%%%%%%%%%%%%%%%%%%%%%%%
%%% BEGIN MATHBF
 \def\bfa{\mathbf a} \def\bfb{\mathbf b} \def\bfc{\mathbf c} \def\bfd{\mathbf d} \def\bfe{\mathbf e} \def\bff{\mathbf f}
 \def\bfg{\mathbf g} \def\bfh{\mathbf h} \def\bfi{\mathbf i} \def\bfj{\mathbf j} \def\bfk{\mathbf k} \def\bfl{\mathbf l}
 \def\bfm{\mathbf m} \def\bfn{\mathbf n} \def\bfo{\mathbf o} \def\bfp{\mathbf p} \def\bfq{\mathbf q} \def\bfr{\mathbf r}
 \def\bfs{\mathbf s} \def\bft{\mathbf t} \def\bfu{\mathbf u} \def\bfv{\mathbf v} \def\bfw{\mathbf w} \def\bfx{\mathbf x}
 \def\bfy{\mathbf y} \def\bfz{\mathbf z} \def\bfA{\mathbf A} \def\bfB{\mathbf B} \def\bfC{\mathbf C} \def\bfD{\mathbf D}
 \def\bfE{\mathbf E} \def\bfF{\mathbf F} \def\bfG{\mathbf G} \def\bfH{\mathbf H} \def\bfI{\mathbf I} \def\bfJ{\mathbf J}
 \def\bfK{\mathbf K} \def\bfL{\mathbf L} \def\bfM{\mathbf M} \def\bfN{\mathbf N} \def\bfO{\mathbf O} \def\bfP{\mathbf P}
 \def\bfQ{\mathbf Q} \def\bfR{\mathbf R} \def\bfS{\mathbf S} \def\bfT{\mathbf T} \def\bfU{\mathbf U} \def\bfV{\mathbf V}
 \def\bfW{\mathbf W} \def\bfX{\mathbf X} \def\bfY{\mathbf Y} \def\bfZ{\mathbf Z} \def\bfw{\mathbf w}
 %%% END MATHBF
%%%%%%%%%%%%%%%%%%%%%%%%%%%%%%% %%%%%%%%%%%%%%%%%%%%%%%%%%%%%%%%%
 %%% BEGIN MATHBB
 \def\R {{\mathbb R }} \def\C {{\mathbb C }} \def\Z{{\mathbb Z}} \def\H{{\mathbb H}} \def\K{{\mathbb K}}
 \def\N{{\mathbb N}} \def\Q{{\mathbb Q}} \def\A{{\mathbb A}} \def\T{\mathbb T} \def\P{\mathbb P} \def\G{\mathbb G}
 \def\bbA{\mathbb A} \def\bbB{\mathbb B} \def\bbD{\mathbb D} \def\bbE{\mathbb E} \def\bbF{\mathbb F} \def\bbG{\mathbb G}
 \def\bbI{\mathbb I} \def\bbJ{\mathbb J} \def\bbL{\mathbb L} \def\bbM{\mathbb M} \def\bbN{\mathbb N} \def\bbO{\mathbb O}
 \def\bbP{\mathbb P} \def\bbQ{\mathbb Q} \def\bbS{\mathbb S} \def\bbT{\mathbb T} \def\bbU{\mathbb U} \def\bbV{\mathbb V}
 \def\bbW{\mathbb W} \def\bbX{\mathbb X} \def\bbY{\mathbb Y} \def\kappa{\varkappa} \def\epsilon{\varepsilon}
 \def\phi{\varphi} \def\le{\leqslant} \def\ge{\geqslant}

\def\UU{\bbU}
\def\Mat{\mathrm{Mat}}
\def\tto{\rightrightarrows}

\def\Gms{\mathrm {Gms}}
\def\Ams{\mathrm {Ams}}
\def\Isom{\mathrm {Isom}}

\def\Gr{\mathrm{Gr}}

\def\graph{\mathrm{graph}}

\def\O{\mathrm{O}}

\def\la{\langle}
\def\ra{\rangle}

%\begin{document}

 \def\ov{\overline}
\def\wt{\widetilde}

\renewcommand{\Re}{\mathop {\mathrm {Re}}\nolimits}
\def\Br{\mathrm {Br}}

 %%%1. ClASSICAL GROUPS
 \def\Isom{\mathrm {Isom}}
\def\SL{\mathrm {SL}}
\def\SU{\mathrm {SU}}
\def\GL{\mathrm {GL}}
\def\U{\mathrm U}
\def\OO{\mathrm O}
 \def\Sp{\mathrm {Sp}}
  \def\GLO{\mathrm {GLO}}
 \def\SO{\mathrm {SO}}
\def\SOS{\mathrm {SO}^*}
 \def\Diff{\mathrm{Diff}}
 \def\Vect{\mathfrak{Vect}}
\def\PGL{\mathrm {PGL}}
\def\PU{\mathrm {PU}}
\def\PSL{\mathrm {PSL}}
\def\Symp{\mathrm{Symp}}
\def\ASymm{\mathrm{Asymm}}
\def\Asymm{\mathrm{Asymm}}
\def\Abs{\mathrm{Abs}}
\def\Gal{\mathrm{Gal}}
\def\End{\mathrm{End}}
\def\Mor{\mathrm{Mor}}
\def\Aut{\mathrm{Aut}}
 \def\PB{\mathrm{PB}}
\def\Ver{\mathrm{Vert}} 
\def\Hier{\mathrm{Hier}} 
\def\Abs{\mathrm{Abs}} 
 
 \def\cA{\mathcal A}
\def\cB{\mathcal B}
\def\cC{\mathcal C}
\def\cD{\mathcal D}
\def\cE{\mathcal E}
\def\cF{\mathcal F}
\def\cG{\mathcal G}
\def\cH{\mathcal H}
\def\cJ{\mathcal J}
\def\cI{\mathcal I}
\def\cK{\mathcal K}
 \def\cL{\mathcal L}
\def\cM{\mathcal M}
\def\cN{\mathcal N}
 \def\cO{\mathcal O}
\def\cP{\mathcal P}
\def\cQ{\mathcal Q}
\def\cR{\mathcal R}
\def\cS{\mathcal S}
\def\cT{\mathcal T}
\def\cU{\mathcal U}
\def\cV{\mathcal V}
 \def\cW{\mathcal W}
\def\cX{\mathcal X}
 \def\cY{\mathcal Y}
 \def\cZ{\mathcal Z}
%%% END MATHCAL %%%%%%%%%%%%%%%%%%%%%%%%%%%%%%%%% %%%%%%%%%%%%%%%%%%%%%%%%%%%%%%%% %%%
\def\0{{\ov 0}}
 \def\1{{\ov 1}}
 %%%%%%%%%%%%%%%%%%%%%%%%%%%% %%%%%%%%%%%%%%%%%%%%%%%%%%%%%%%%%%% %%% BEGIN GOTIC
 \def\frA{\mathfrak A}
 \def\frB{\mathfrak B}
\def\frC{\mathfrak C}
\def\frD{\mathfrak D}
\def\frE{\mathfrak E}
\def\frF{\mathfrak F}
\def\frG{\mathfrak G}
\def\frH{\mathfrak H}
\def\frI{\mathfrak I}
 \def\frJ{\mathfrak J}
 \def\frK{\mathfrak K}
 \def\frL{\mathfrak L}
\def\frM{\mathfrak M}
 \def\frN{\mathfrak N} \def\frO{\mathfrak O} \def\frP{\mathfrak P} \def\frQ{\mathfrak Q} \def\frR{\mathfrak R}
 \def\frS{\mathfrak S} \def\frT{\mathfrak T} \def\frU{\mathfrak U} \def\frV{\mathfrak V} \def\frW{\mathfrak W}
 \def\frX{\mathfrak X} \def\frY{\mathfrak Y} \def\frZ{\mathfrak Z} \def\fra{\mathfrak a} \def\frb{\mathfrak b}
 \def\frc{\mathfrak c} \def\frd{\mathfrak d} \def\fre{\mathfrak e} \def\frf{\mathfrak f} \def\frg{\mathfrak g}
 \def\frh{\mathfrak h} \def\fri{\mathfrak i} \def\frj{\mathfrak j} \def\frk{\mathfrak k} \def\frl{\mathfrak l}
 \def\frm{\mathfrak m} \def\frn{\mathfrak n} \def\fro{\mathfrak o} \def\frp{\mathfrak p} \def\frq{\mathfrak q}
 \def\frr{\mathfrak r} \def\frs{\mathfrak s} \def\frt{\mathfrak t} \def\fru{\mathfrak u} \def\frv{\mathfrak v}
 \def\frw{\mathfrak w} \def\frx{\mathfrak x} \def\fry{\mathfrak y} \def\frz{\mathfrak z} \def\frsp{\mathfrak{sp}}
 %% This is Lie algebra %%% END GOTIC
%%%%%%%%%%%%%%%%%%%%%%%%%%%%%%%% %%%%%%%%%%%%%%%%%%%%%%%%%%%%%%%%%
%%% BEGIN MATHBF
 \def\bfa{\mathbf a} \def\bfb{\mathbf b} \def\bfc{\mathbf c} \def\bfd{\mathbf d} \def\bfe{\mathbf e} \def\bff{\mathbf f}
 \def\bfg{\mathbf g} \def\bfh{\mathbf h} \def\bfi{\mathbf i} \def\bfj{\mathbf j} \def\bfk{\mathbf k} \def\bfl{\mathbf l}
 \def\bfm{\mathbf m} \def\bfn{\mathbf n} \def\bfo{\mathbf o} \def\bfp{\mathbf p} \def\bfq{\mathbf q} \def\bfr{\mathbf r}
 \def\bfs{\mathbf s} \def\bft{\mathbf t} \def\bfu{\mathbf u} \def\bfv{\mathbf v} \def\bfw{\mathbf w} \def\bfx{\mathbf x}
 \def\bfy{\mathbf y} \def\bfz{\mathbf z} \def\bfA{\mathbf A} \def\bfB{\mathbf B} \def\bfC{\mathbf C} \def\bfD{\mathbf D}
 \def\bfE{\mathbf E} \def\bfF{\mathbf F} \def\bfG{\mathbf G} \def\bfH{\mathbf H} \def\bfI{\mathbf I} \def\bfJ{\mathbf J}
 \def\bfK{\mathbf K} \def\bfL{\mathbf L} \def\bfM{\mathbf M} \def\bfN{\mathbf N} \def\bfO{\mathbf O} \def\bfP{\mathbf P}
 \def\bfQ{\mathbf Q} \def\bfR{\mathbf R} \def\bfS{\mathbf S} \def\bfT{\mathbf T} \def\bfU{\mathbf U} \def\bfV{\mathbf V}
 \def\bfW{\mathbf W} \def\bfX{\mathbf X} \def\bfY{\mathbf Y} \def\bfZ{\mathbf Z} \def\bfw{\mathbf w}
 %%% END MATHBF
%%%%%%%%%%%%%%%%%%%%%%%%%%%%%%% %%%%%%%%%%%%%%%%%%%%%%%%%%%%%%%%%
 %%% BEGIN MATHBB
 \def\R {{\mathbb R }} \def\C {{\mathbb C }} \def\Z{{\mathbb Z}} \def\H{{\mathbb H}} \def\K{{\mathbb K}}
 \def\N{{\mathbb N}} \def\Q{{\mathbb Q}} \def\A{{\mathbb A}} \def\T{\mathbb T} \def\P{\mathbb P} \def\G{\mathbb G}
 \def\bbA{\mathbb A} \def\bbB{\mathbb B} \def\bbD{\mathbb D} \def\bbE{\mathbb E} \def\bbF{\mathbb F} \def\bbG{\mathbb G}
 \def\bbI{\mathbb I} \def\bbJ{\mathbb J} \def\bbL{\mathbb L} \def\bbM{\mathbb M} \def\bbN{\mathbb N} \def\bbO{\mathbb O}
 \def\bbP{\mathbb P} \def\bbQ{\mathbb Q} \def\bbS{\mathbb S} \def\bbT{\mathbb T} \def\bbU{\mathbb U} \def\bbV{\mathbb V}
 \def\bbW{\mathbb W} \def\bbX{\mathbb X} \def\bbY{\mathbb Y} \def\kappa{\varkappa} \def\epsilon{\varepsilon}
 \def\phi{\varphi} \def\le{\leqslant} \def\ge{\geqslant}

\def\UU{\bbU}
\def\Mat{\mathrm{Mat}}
\def\tto{\rightrightarrows}

\def\Gr{\mathrm{Gr}}

\def\graph{\mathrm{graph}}

\def\O{\mathrm{O}}

\def\la{\langle}
\def\ra{\rangle}

\begin{center}
\bf\large
The subgroup $\PSL_2(\R)$ is spherical in the group of diffeomorphisms of the circle

\sc Yury A. Neretin%
\footnote{Supported by the grant FWF, Project P25142}
\end{center}

{\small We show that the group $\PSL_2(\R)$ is a spherical subgroup
in the group of $C^3$-diffeomorphisms
of the circle. Also, the group of automorphisms of a Bruhat--Tits tree
is a spherical subgroup in the group of hierarchomorphisms of the tree.}

\sm

{\bf \punct Sphericity.} Let $G$ be a topological group, $K$ be a subgroup.
An irreducible unitary representation $\rho$ of $G$ in a Hilbert space
$H$ is called {\it spherical}
if there is a unique up to a scalar factor non-zero $K$-invariant vector $v$ in $H$.
The matrix element 
$$
\Phi(g):=\la\rho(g)v, v\ra
$$
is called a {\it spherical function} on $G$.
A subgroup $K$ in $G$ is called {\it spherical} if for any irreducible unitary representation
of $G$ the dimension of the space of $K$-invariant vectors is $\le 1$.

\sm

For  various types of spherical pairs in this sense, see \cite{Gel},
\cite{Kra}, \cite{Olsh-GB}, \cite{Cec}, \cite{Ness}, \cite{Ner-spheric}. For all known examples
the group $K$ is compact or is  an infinite-dimensional analog of compact 
groups as $\U(\infty)$, $\O(\infty)$, $\Sp(\infty)$, $\mathrm S(\infty)$ etc.
('heavy groups' in the sense of \cite{Ner-book}).

\sm

{\bf \punct Statements.} Let $\SL_2(\R)$ be the group of $2\times 2$
real matrices with determinant 1, let $\PSL_2(\R)$ be its quotient with respect
to the center, $\SL_2(\R)^\sim$ be the universal covering group.
Denote by $\Diff$ (respectively by $\Diff^{3}$) the group of $C^\infty$-smooth (resp. $C^3$-smooth) 
orientation preserving diffeomorphisms
of the circle. Denote by $\Diff^\sim$ the universal covering of $\Diff$, we realize
$\Diff^\sim$ as the group of smooth diffeomorphisms $q$ of the line $\R$
satisfying the condition
$$
q(\phi+2\pi)=q(\phi)+2\pi.
$$
The {\it Bott cocycle} $c(q_1,q_2)$ on  $\Diff^\sim$ is defined
by the formula
$$
c(q_1,q_2)=\int_0^{2\pi}\ln q_1'(q_2(\phi) d\ln q_2'(\phi).
$$
Consider the central  extension $\wt \Diff$  of $\Diff^\sim$
determined by the Bott cocycle (see, e.g., \cite{Fuks}, \S3.4).
 By $\wt\Diff^{3}$ we denote the similar central extension
of $\Diff^{3}$.

\begin{theorem}
The subgroup $\PSL_2(\R)$ is spherical in the group $\Diff^3$.
\end{theorem}

\begin{theorem}
The subgroup $\PSL_2^\sim(\R)$ is spherical in the group $\wt\Diff^3$.
\end{theorem}

{\sc Remark.} a) Several series of spherical representations of $\Diff$ 
and $\Diff^\sim$ 
were constructed in \cite{Ner-dan}, see also \cite{Ner-book}, \S IX.6.
All such spherical representations are continuous in the $C^3$-topology.

\sm

b)  Sobolev diffeomorphisms of the circle of the class $s>3/2$
form a group (see \cite{IK}, Theorem 1.2 and Appendix B). Our proof survives
for the group of Sobolev diffeomorphisms of the class $H^{4.5}$.
\hfill $\boxtimes$.

\sm

Next, consider a combinatorial analog of $\Diff$.
Fix an integer $n\ge2 $. Consider the {\it Bruhat--Tits tree} $T_n$, i.e., the infinite tree 
such that each vertex is incident to $n+1$ edges. Let $\Abs(T_n)$ be its boundary
(for detailed definitions, see, e.g., \cite{Olsh}, \cite{Ner-tree}). Denote
by $\Aut(T_n)$ the group of all automorphisms of the graph $T_n$.
It is a locally compact group, stabilizers of finite subtrees form a base
of open-closed neighborhoods of unit.

Denote by $\Ver(T_n)$ the set of vertices of $\T_n$. Consider a bijection
$\theta:\Ver(T_n)\to\Ver(T_n)$
such that for all but a finite numbers of pairs of adjacent vertices $(a,b)$,
vertices $\theta(a)$, $\theta(b)$ are adjacent.
  {\it Hierarchomorphism} of the tree $T_n$ is are homeomorphisms of $\Abs(T_n)$
induced by such  maps, see \cite{Ner-tree}, \cite{Ner-h}. Denote by
$\Hier(T_n)$ the group of all hierarchomorphisms of the tree $T_n$.

\sm

{\sc Remark.} a) For a prime $n=p$ the boundary $\Abs(T_p)$ can be identified 
with a $p$-adic projective line. The group $\Aut(T_p)$ contains the $p$-adic
$\PSL_2$ and the representation theory of $\Aut(T_p)$ is similar to 
the representation theory
of $p$-adic and real $\SL_2$ (see \cite{Car}, \cite{Olsh}). 
The group $\Hier(T_p)$  contains  the group of locally analytic
diffeomorphisms of the $p$-adic projective line.

b) Richard Thompson groups (see \cite{RT}) are discrete subgroups of  $\Hier(T_n)$.

c) Several series of $\Aut(T_n)$-spherical representations
of $\Hier(T_n)$ were constructed in \cite{Ner-p}, \cite{Ner-tree}.
\hfill $\boxtimes$

\sm

We define a topology on the group $\Hier(T_n)$ assuming that $\Aut(T_n)$ is an open subgroup
(the coset space $\Hier(T_n)/\Aut(T_n)$ is countable).

\begin{theorem}
The subgroup $\Aut(T_p)$ is spherical in $\Hier(T_p)$.
\end{theorem}

%The following theorem is semi-trivial.

\begin{theorem}
Let $G\supset K$ be a spherical pair. Assume that $K$ does not admit
nontrivial finite-dimensional unitary representations.
 Let $\Phi_1(g)$, $\Phi_2(g)$ be $K$-spherical functions on $G$. Then
$\Phi_1(g)\Phi_2(g)$ is a spherical function.
\end{theorem}

The both $K=\SL_2(\R)^\sim$ and $\Aut(T_n)$ satisfy this condition.

\sm

{\bf \punct Proof of Theorem 1.}
Fix a point $a$ in the circle. Denote by $G_0\subset \Diff$ the group of diffeomorphsims
$q$ such that $q(x)=x$ in a neighborhood of $a$. By $G^*$ we denote the group
of diffeomorphisms that are flat at $a$, i.e., 
$$q(a)=a,\quad q'(a)=1,\quad q''(a)=q'''(a)=\dots=0$$

Let $\rho$ be an irreducible unitary representation of $\Diff$ in $H$. Denote by
$V$ the subspace of all $\PSL_2(\R)$-fixed vectors. Let $P$ be the operator of orthogonal
projection on $V$. For $h\in \PSL_2(\R)$ we have
$$
P\rho(h)=\rho(h)P=P.
$$  
Denote
$$
\wh\rho(g):=P\rho(g)P.
$$
This operator depends only on a double coset
of $\Diff$ by $\PSL_2(\R)$,
$$
\wh\rho(h_1gh_2):=\wh\rho(g), \qquad h_1,h_2\in\PSL_2(\R).
$$

\begin{lemma}
If $\rho$ is continuous in the $C^3$-topology, then the
operators $\wh \rho(g)$ pairwise commute.
\end{lemma}

{\sc Proof.}
The following statement is our key argument: 

{\it Let a sequence $h_j\in \PSL_2(\R)$ converges to infinity%
\footnote{I.e., for any compact subset $B$, we have $h_j\notin B$ starting some number.}. 
 Then $\rho(h_j)$ weakly converges to $P$,} see Howe, Moore \cite{HM}, Theorem 5.1
(this is a general theorem for semisimple group,
for $\PSL_2(\R)$ it can be easily verified case-by-case).

Let us realize the circle as the real projective line $\R\bbP^1=\R\cup\infty$. Without loss
of generality we can set $a=\infty$. Let $U_t(x)=x+t$ be a shift on $\R$, we have
$U_t\in \PSL_2(\R)$. Consider diffeomorphisms $r$, $q\in G_0$. For sufficiently large
$t$ the supports of $r$ and $U_t\circ q\circ U_{-t}$ are disjoint.
Therefore, these diffeomorphisms commute. Hence,
$$
\rho(r)\,\rho(U_t)\,\rho(q)\,\rho(U_{-t})=\rho(U_t)\,\rho(q)\,\rho(U_{-t})\,\rho(r)
.
$$
Therefore,
$$
P\,\rho(r)\,\rho(U_t)\,\rho(q)\,P=P\,\rho(q)\,\rho(U_{-t})\,\rho(r)\,P
.
$$
Passing to a weak limit as $t\to+\infty$, we get
$$
P\,\rho(r)\,P\,\rho(q)\,P=P\,\rho(q)\,P\,\rho(r)\,P
.
$$
Thus
\begin{equation*}
\wh \rho(r)\,\wh \rho(q)= \wh \rho(q)\,\wh \rho(r), \qquad
\text{where $r$, $q\in G_0$.}
\end{equation*}
But $G_0$ is dense in $G^*$. Therefore the same identity holds for $r$, $q\in G_*$.
Indeed, let $r_j$, $q_j\in G_0$ be sequences convergent to $r$, $q$ respectively.
Passing to the iterated limit
$$
\lim_{j\to\infty}\Bigl(\lim_{k\to\infty} \wh \rho(r_j)\,\wh \rho(q_k)\Bigr)=
\lim_{j\to\infty}\Bigl(\lim_{k\to\infty} \wh \rho(q_k) \wh \rho(r_j)\Bigr)
$$
and keeping in mind the separate weak continuity of the operator product,
we get the desired statement.

Our last argument: {\it the set $\PSL_2(\R)\cdot G_* \cdot \PSL_2(\R)$
is dense in $\Diff$ with respect to the $C^3$-topology.}

Let us prove this.
Choose a coordinate on $\R\bbP^1$
 such that $a=0$. 
Let $q\in \Diff$. Consider its Schwarzian derivative,
$$
S(q)=\frac{q'q'''-\frac32 (q'')^2}{(q')^2}
.
$$
Consider a point $b$ such that $S(q)(b)=0$ 
(by the Ghys theorem, the Schwarzian derivative of a diffeomorphisms of the circle
 has at least 4 zero, see \cite{OT}, Theorem 4.2.1). Then for 
 $$r:=U_{-q(b)}\circ q \circ U_{b}$$
 we have
 $ r(0)=0$, $S(r)(0)=0$. Consider  maps 
 $$
 \sigma(x)=\frac{ux}{u^{-1}+vx}
 ,$$ 
such $\sigma\in \PSL_2(\R)$ fix 0.
Choosing parameters $u$, $v$, we can achieve 
$$(r\circ \sigma)'(0)=1,\qquad (r\circ \sigma)''(0)=0.$$
 Recall the transformation  property of the Scwarzian:
$$
S(r\circ\sigma)=(S(r)\circ\sigma)\cdot (\sigma')^2+ S(\sigma).
$$
Since $\sigma$ is linear fractional, $S(\sigma)=0$.
Therefore $S(r\circ\sigma)=0$, and $(r\circ\sigma)'''(0)=0$.
Such $r\sigma$ can be approximated in $C^3$-topology by elements of $G_*$.
This proves Lemma 5. \hfill $\square$

\sm

Theorem 1 is a corollary of the lemma. Note, that
$\wh \rho(g)^*=\wh \rho(g^{-1})$. Thus we get a family of commuting operators
in $V$, such that an adjoint operator $A^*$ is contained in the family
together with $A$. If $\dim V>1$, then this family has a proper invariant
 subspace in $V$, say $W$. Consider the $\Diff$-cyclic span of $W$, i.e., 
 the subspace $Z$ spanned by vectors $\rho(g) w$, where $g\in\Diff^3$ and
 $w\in W$. Then
 $$P\rho(g)w=P\rho(g)P\,w=\wh\rho(g)w\in W.$$
  Hence,  $PZ=W$ and therefore $Z$   is a proper subspace in $H$.

  \sm

{\bf \punct Proof of Theorem 2.}
It  repeats the previous proof with two additional remarks. 

1) Consider the homomorphism 
$\pi:\SL_2(\R)^\sim\to \PSL_2(\R)\simeq \SL_2(\R)^\sim/\Z$.
 We say that a sequence $h_j\in  \SL_2(\R)^\sim$ converges to $\infty$ if 
 $\pi (h_j)\to\infty$. Then the Howe--More theorem remains  valid.
 
2) For a pair of diffeomorphisms
 with disjoint supports $p$, $q$ the Bott cocycle $c(q,p)$ vanishes,
 hence the diffeomorphisms $p$, $q$ commute in the extended 
 group.

 \sm
 
  {\bf \punct Proof of Theorem 3.} First, there is the following analog of the Howe--Moore
  theorem:
{\it Let a sequence $h_j\in\Aut(T_n)$ converges to $\infty$. Then for any unitary
representation $\rho$ of $\Aut(T_n)$ the sequence $\rho(h_j)$ converges to the projection
to the subspace of $\Aut(T_n)$-fixed vectors,} see \cite{Lub};  this can be easily
verified case-by-case starting the classification theorem of \cite{Olsh}.

  Second, fix a point $a\in\Abs(T_n)$ and denote by $G_0$ the group of hierarchomorphisms
  that are trivial in a neighborhood of $a$. Let $q$, $r\in G_0$. Then there is 
  a sequence $h_j\in \Aut(T_n)\cap G_0$ such that $h_j$ tends to $\infty$ and supports 
  of $h_j p h_j^{-1}$ and $q$ are disjoint. We omit a proof, since it is easier 
  to understand its self-evidence than to read a formal exposition.
   
   Third, 
   $$
   \Aut(T_n)\cdot G_0 \cdot \Aut(T_n)=\Hier(T_p)
   $$
   
Now we can repeat the proof of Theorem 1.   
 
\sm 
 
 {\bf \punct Proof of Theorem 4.} The statement is semitrivial.
 
 \begin{lemma}
Let $\nu_1$ $\nu_2$ be unitary representations of a group $\Gamma$. If
the tensor product $\nu_1\otimes\nu_2$ contains a nonzero $\Gamma$-invariant vector,
then the both $\nu_1$ and $\nu_2$ have finite-dimensional subrepresentations.  
 \end{lemma}
 
{\sc Proof of the lemma.} Assume that an invariant vector exists.
 Denote the spaces of representations by $V_1$, $V_2$.
We identify $V_1\otimes V_2$ with the space of Hilbert--Schmidt operators
$V_1'\to V_2$, where $V_1'$ is the dual space to $V_1$. An invariant vector 
corresponds to an intertwining operator $T:V_1'\to V_2$.
The operator $TT^*$ is an intertwining operator in $V_2$. Since $TT^*$
is compact and nonzero, 
it has a finite-dimensional  eigenspace, and this subspace is $G$-invariant.
\hfill $\square$

 \sm
 
{\sc Proof of the theorem}.
 Let $\rho_1$ and $\rho_2$ be $K$-spherical representations of $G$ in $H_1$ and $H_2$.
 Let $v_1$, $v_2$ be fixed vectors. By the lemma, $v_1\otimes v_2$ is a unique
 $K$-fixed vector in $H_1\otimes H_2$. The cyclic span
 $W$ of $v_1\otimes v_2$ is an irreducible subreprepresentation.
 Indeed, let $W=W_1\oplus W_2$ be
 reducible. Then the both projections of $v_1\otimes v_2$ to $W_1$, $W_2$
 are $K$-fixed, therefore $v_1\otimes v_2$ must be contained
in one of summands, say $W_1$, and thus the cyclic span of $v_1\otimes v_2$
is contained in $W_1$, i.e., $W=W_1$.

Now we consider the representation of $G$ in $W$,
\begin{multline*}
\left\la\bigl(\rho_1(g)\otimes \rho_2(g)\bigr) v_1\otimes v_2, v_1\otimes v_2\right\ra_W
=\la \rho_1(g)v_1,v_1\ra_{H_1}\cdot\la \rho_2(g)v_2,v_2\ra_{H_2}
\\=
\Phi_1(g)\Phi_2(g).
\end{multline*}

\noindent
\tt Math.Dept., University of Vienna,
 \\
 Oskar-Morgenstern-Platz 1, 1090 Wien;
 \\
\& Institute for Theoretical and Experimental Physics (Moscow);
\\
\& Mech.Math.Dept., Moscow State University.
\\
e-mail: neretin(at) mccme.ru
\\
URL:www.mat.univie.ac.at/$\sim$neretin

\end{document}